\documentclass[11pt,a4paper]{article}
\usepackage[utf8]{inputenc}
\usepackage{amsmath}
\usepackage{amsfonts}
\usepackage{amssymb}
\usepackage{bm}
\usepackage{pdflscape}
\usepackage[margin=1in,footskip=0.25in]{geometry}
\usepackage{graphicx}
\usepackage{multirow}
\usepackage{graphicx}
\usepackage{wrapfig}
\usepackage{setspace}
\usepackage{lscape}
\usepackage[table,xcdraw]{xcolor}
\usepackage{float}

\usepackage{longtable}
\usepackage{courier}
\usepackage[printwatermark]{xwatermark}
\begin{document}

\noindent\textbf{\Large{Modified estimator for the proportion of true null hypotheses under discrete setup with proven FDR control by the adaptive Benjamini-Hochberg procedure}}\\

\begin{center}
Aniket Biswas\\
Department of Statistics\\
Dibrugarh University\\
Dibrugarh, Assam, India-786004\\
Email:\textit{biswasaniket44@gmail.com}
\end{center}

\begin{center}
Gaurangadeb Chattopadhyay\\
Department of Statistics\\
University of Calcutta\\
Kolkata, West Bengal, India-700019\\
Email:\textit{gcdhstat@gmail.com}
\end{center}

\section*{Abstract}
Some crucial issues about a recently proposed estimator for the proportion of true null hypotheses ($\pi_0$) under discrete setup are discussed. An estimator for $\pi_0$ is introduced under the same setup. The estimator may be seen as a modification of a very popular estimator for $\pi_0$, originally proposed under the assumption of continuous test statistics. It is shown that adaptive Benjamini-Hochberg procedure remains conservative with the new estimator for $\pi_0$ being plugged in.\\

\noindent \textbf{Keywords}: Multiple hypotheses testing, RNA Sequence data, Fisher's exact test.\\

\noindent \textbf{MSC 2010:} 62F03, 62P10.

\section{Introduction}
In this era of next generation sequencing data, application of multiple testing procedures in discrete paradigm have become increasingly popular over the last decade. The classical Benjamini-Hochberg (henceforth BH) procedure controls FDR at a prefixed level when the $p$-values corresponding to the true null hypotheses are uniform (Benjamini and Hochberg 1995). The additional assumptions for proving the control over FDR are that the $p$-values corresponding to the true null hypotheses are independent among themselves and they are independent with the $p$-values corresponding to the false null hypotheses. Following the steps in the proof of FDR control by BH procedure, it is worth noting that the procedure is conservative in two dimensions. The first one is that the BH procedure originally controls FDR at level $\pi_0 q$ instead of the prefixed level $q\in (0,1)$ and thus a reasonable estimate of $\pi_0$ can inflate the cutoffs for each of the ordered $p$-values to bring more rejections with the same control over FDR. It is worth mentioning that, the BH procedure was originally introduced for multiple continuous test-statistics where the tested nulls are usually simple. Despite the above fact, BH procedure accommodates multiple discrete test statistics as the $p$-values under this setup are stochastically larger than the uniform variate. Unfortunately, super-uniformity of the $p$-values corresponding to the true null hypotheses makes the procedure conservative under discrete setup. Heyse (2011) noted the fact and modified the BH procedure (henceforth BHH procedure) for discrete test statistics. However, the first cause of being conservative is valid for both continuous and discrete paradigms. A huge literature is available on estimating $\pi_0$ under continuous set-up (Storey 2002, Storey et al. 2004, Langaas et al. 2005, Benjamini et al. 2006, Wang et al. 2011, Cheng et al. 2015, Biswas (1) 2020 among many others). These estimators for $\pi_0$ become very conservative when applied in discrete setup. Dialsingh et al. (2015) and Chen et al. (2018) are the only works particularly focusing on the issue of estimating $\pi_0$ under discrete setup till date. Suitable estimate of $\pi_0$ is used to construct more powerful adaptive FDR controlling procedures (Benjamini et al. 2006, Sarkar 2008, Blanchard and Roquain 2009, Chen and Doerge 2014). The current work introduces a new estimator for the proportion of true null hypotheses under discrete setup. The focus is obviously on the desired result that the adaptive FDR controlling algorithms constructed by plugging in the new estimator for $\pi_0$ is conservative. The work of Chen et al. (2018) and two related reports need a revisit to motivate the current work.\\

The rest of the article is organized as follows. In the next section we introduce some useful notations, discuss some key developments and motivate the current research. We present the estimator in Section 3. The theoretical result on adaptive algorithms being conservative is stated and proved in Section 4. We conclude the article by pointing out the possible future additions to the article.  

\section{Background}
Suppose $H_1, H_2, ..., H_m$ are $m$ related but independent hypotheses to be tested. From empirical Bayesian motivation given in Storey (2002), $H_i$'s are iid Bernoulli random variables for $i\in \mathcal{I}=\{1, 2, ..., m\}$ with parameter $\pi_0$ and $H_i=1$ denotes that the $i$-th null hypothesis is true. Let $\mathcal{T}=\{i\in \mathcal{I}: H_i=1\}$, the set of indices corresponding to the true null hypotheses and $\mathcal{F}=\{i \in \mathcal{I}: H_i=0\}$, the same corresponding to the false null hypotheses. We denote the number of true null hypotheses by $m_0$ and $m_0=\sum_{i \in \mathcal{I}}H_i\sim$ Binomial$(m,\pi_0)$. $H_i$'s and hence $m_0$ are unobserved. The parameter of interest $\pi_0$ is required and it is usually estimated from the set of available $p$-values $\{p_1, p_2, ..., p_m\}$, $p_i$ being the $p$-value corresponding to $H_i$, $i\in \mathcal{I}$. We use the notation $p$ to denote $p$-value irrespective of it being a random variable or an observed value and the meaning holds according to the context. Under discrete multiple testing setup, multiple Fisher's exact tests (FET) or Binomial tests (BT) are performed. For details on two-population FET, BT and their relevance in identifying differentially expressed genes from RNA sequence data available on two different study groups, see Chen (2019). \\

For a suitably chosen tuning parameter $\tau\in (0,1)$, Storey (2002) assumed that $\{i\in \mathcal{F}: p_i>\tau\}$ has cardinality $0$. Let $I$ denote the usual indicator function such that $I(x)=1$ if the statement $x$ is true and $I(x)=0$, otherwise. As discussed earlier, $p\sim$ Uniform$(0,1)$ when the tested null is simple and the test statistics are continuous. Thus under continuous setup with the additional assumption, $\sum_{i\in \mathcal{I}} I(p_i>\tau)\sim$ Binomial$(m_0, 1-\tau)$. Noting that $E(m_0)=m\pi_0$, a conservative estimator for $\pi_0$ is 
\begin{equation}\label{eq:pi0_Storey2002}
    \hat{\pi}_0(\tau)=\frac{1}{m(1-\tau)} \sum_{i\in \mathcal{I}} I(p_i>\tau).
\end{equation}
Storey et al. (2004) further introduced a conservative bias in $\hat{\pi}_0 (\tau)$ and the modified estimator is 
\begin{equation}\label{eq:pi0_Storey2004}
    \hat{\pi}_0^S(\tau)=\frac{1}{m(1-\tau)}+\frac{1}{m(1-\tau)} \sum_{i\in \mathcal{I}} I(p_i>\tau).
\end{equation}
Though the distributions of the true null $p$-values are same under the continuous setup, distribution of the false null $p$-values are different and thus should be labelled by $i$. Denote the distribution function of $p_i$ for $i\in\mathcal{F}$ by $G_i$. Then $\hat{\pi}_0^S(\tau)$ has conservative bias
\begin{equation}\label{eq:Storeybias_continuous}
    B_1=\frac{1}{m(1-\tau)}+\frac{1}{m(1-\tau)}\sum_{i\in \mathcal{F}}[1-G_i(\tau)].
\end{equation}
However, if $\hat{\pi}_0^S(\tau)$ is used for estimating $\pi_0$ under discrete setup, another positive bias term gets introduced. For discrete tests, the true null $p$-values lose homogeneity and we denote distribution function of $p_i$ for $i \in \mathcal{T}$ by $F_i$. Under discrete multiple testing setup, the bias in $\hat{\pi}_0^S(\tau)$ is
\begin{equation}\label{eq:Storeybias_discrete}
B_2=\frac{1}{m(1-\tau)}+\frac{1}{m(1-\tau)}\sum_{i\in \mathcal{T}}[\tau-F_i(\tau)] +\frac{1}{m(1-\tau)}\sum_{i\in \mathcal{F}}[1-G_i(\tau)].
\end{equation}
Henceforth we continue our discussion by assuming that the $p$-values are obtained from discrete tests. As mentioned in Section 1, $p$ stochastically dominates Uniform$(0,1)$ random variable when it is obtained from a discrete test. Thus the additional bias $B_2-B_1$ is non-negative since $F_i(\tau)\leq \tau$. Let $S_i$ denote the support of $p_i$. It is to be noted that $B_2=B_1$ iff $\tau\in S_i$ for all $i\in \mathcal{I}$, that is $\tau\in \cap_{i\in \mathcal{I}} S_i$. Storey (2002) suggested a bootstrap routine to select appropriate $\tau$ under continuous setup and Chen et al. (2018) worked with $\tau=0.5$ under discrete setup. Whatever be the technique we follow, the choice of $\tau$ is prefixed and there is no guarantee whatsoever that $\tau$ can be chosen from $\cap_{i\in \mathcal{I}}S_i$. We illustrate the situation with a numerical example for FET using the notations in section 2.1 of Chen (2019). Consider $X_i\sim$ Binomial$(q_i,N)$ for $i=1,2$, with the following testing problem 
\begin{equation}\label{eq:FET}
    H_0: q_1=q_2\quad \textrm{versus}\quad H_1: q_1\neq q_2.
\end{equation}
for discrete tests, two-sided $p$-values are defined in many ways but to maintain super-uniformity of the $p$-values, we work with the usual one  similar to Chen et al. (2018). A concise discussion on super-uniform and sub-uniform $p$-values for discrete test interested readers may refer to Chen (2019). If the total observed count $c=1$, least possible $p$-values is $1$. Thus if a single test with $c=1$, is to be performed among the $m$ tests, $\cap_{i\in \mathcal{I}} S_i$ is essentially $1$. Thus, the only possible choice of $\tau$ for $B_1=B_2$ is $1$ which is clearly an absurd choice. However, for all practical purposes the rows with total observed count $1$ are removed from the data matrix for further analysis. Henceforth we assume that the working dataset does not contain any row with total count equal to $1$. Still the problem persists as follows. Assume $m=3$ FETs, $N=5$ and total observed counts are $2, 3$ and $4$, respectively. Then $S_1=\{0.4444, 1.0000\}$, $S_2=\{0.1667, 1.000\}$ and $S_3=\{0.0476, 0.5238, 1.0000\}$. Again the only possible choice of $\tau$ for $B_1=B_2$ is $1$. Thus it is clear that while testing $m$ hypotheses simultaneously under discrete setup, use of $\hat{\pi}_0^S(\tau)$ amounts to $B_2-B_1$ extra bias eventually.\\

Chen et al. (2018) came up with an idea to not work with prefixed $\tau$. In fact, the authors suggested to use different guiding values for different $p$-values in $\hat{\pi}_0^S(\tau)$, to overcome the above mentioned shortcoming. The idea is to fix tuning parameter $\tau$ and determine a guiding value $\lambda_i\in S_i$ only for $p_i$ depending on the fixed $\tau$. The estimator $\hat{\pi}_0^G$ in Chen et al. (2018) can be  perceived to be an improvement over $\hat{\pi}_0^S$ since use of customized guiding values vanishes the extra bias term as demonstrated therein. Chen et al. (2018) used multiple tuning parameters for possible reduction in variance of the estimator in the same spirit of Jiang and Doerge (2008). The estimator was proven to be able to construct conservative adaptive BH algorithm. Unfortunately, the proof of the result suffers from a critical mistake as pointed out in Biswas (2) (2020). However, a revised version of the proof makes the assumption that the set of tuning parameters $\{\tau_j:j=1, 2, ..., n\}$ are chosen such that $\cup_{j=1}^n\{\tau_j\}\subset \cap_{i\in \mathcal{I}}S_i$ (Chen and Doerge 2020). The authors also mentioned that such a choice may not exist for multiple tuning parameters and thus only option is to work with a single tuning parameter. Under this restrictive assumption on $\tau$ the tuning parameter, $\lambda_i$ as defined in Chen et al. (2018) is identically $\tau\in \cap_{i \in\mathcal{I}}S_i$ for all $i\in \mathcal{I}$. Thus under the assumption $\hat{\pi}_0^G$ becomes identical to $\hat{\pi}_0^S$ for a fixed $\tau$. It is obvious that the shortcomings of $\hat{\pi}_0^S$ discussed in the previous paragraph are also present for the recently developed $\hat{\pi}_0^G$, under the assumption made in the revised proof. The above facts necessitate further investigation. In this article, we propose a new estimator for $\pi_0$ in similar spirit of the work in Chen et al. (2018). Choice of the tuning parameters $\{\tau_j:j=1, 2, ..., n\}$ are not restricted by the support of the $p$-values for the new estimator. However, for each $p_i$ and $\tau_j$, a threshold value $\lambda_{ij}\in S_i$ is obtained in a different way from $\lambda_{ij}$'s in Chen et al. (2018). Despite the liberal choices of tuning parameters we have been able to prove that conservative adaptive FDR controlling step-up procedures can be constructed by plugging in the new estimator. We also point out that under the assumption in Chen and Doerge (2020), the new estimator is same as $\hat{\pi}_0^G$ and hence $\hat{\pi}_0^S$. Though the new estimator is developed from a novel perspective in this article, it may be thought of as a generalization of the estimator in Chen et al. (2018) or a modification of the estimator in Storey et al. (2004) under discrete paradigm with proven FDR control for the adaptive BH procedure.            
\section{Method of estimation}
 Define $q_i=\inf S_i$ for $i\in \mathcal{I}$ and $\nu=\max\{q_i:1\leq i\leq m\}$. We assume that $\nu<1$ as the rows of the dataset with total observed count $c=1$ are to be removed as already mentioned in Section 2. The estimator and the desired result do not require this assumption but we continue with the reasonable assumption for two reasons. The first one is that the data reduction technique makes further algebraic treatment straight forward and the second one is that the technique reduces the number of tests to perform by screening out those rows of the dataset for which testing the corresponding hypothesis is futile due to extremely less available information. One may refer to the data analysis sections of Chen et al. (2018) and Chen (2019) for verification of the discussed usual practice. Now we discuss the estimator through the following steps.
 \subsection*{Algorithm 1}
\begin{itemize}
    \item Set a sequence of tuning parameters $\nu\leq\tau_1\leq\tau_2\leq ... \leq \tau_n<1$. 
    
    \item For each $i\in \mathcal{I}$ and $j\in \mathcal{J}=\{1, 2, ..., n\}$, define
    \begin{eqnarray*}
    & T_{ij}=\{\lambda\in S_i: \lambda\geq \tau_j\}\\
    & \lambda_{ij}=\inf \{\lambda:\lambda\in T_{ij}\}\\
    & \eta_j= \max_{1\leq i\leq m} \lambda_{ij}.
    \end{eqnarray*}
    
    \item For each $j\in \mathcal{J}$, define the following trial estimator for $\pi_0$.
    \begin{equation*}
        \beta(\tau_j)=\frac{1}{m(1-\eta_j)}+\frac{1-\tau_j}{m(1-\eta_j)}\sum_{i\in \mathcal{I}}\frac{I(p_i>\lambda_{ij})}{1-\lambda_{ij}}
    \end{equation*}
    If $\beta(\tau_j)$ is greater than $1$, we take $\beta(\tau_j)=1$.
    \item The final estimator for $\pi_0$ is 
    \begin{equation*}
        \hat{\pi}_0^H=\frac{1}{n}\sum_{j\in \mathcal{J}}\beta(\tau_j).
    \end{equation*}
\end{itemize}

\noindent Note that the assumption $\cup_{i\in \mathcal{J}}\{\tau_j\}\subset \cap_{i \in \mathcal{I}}S_i\implies\tau_j=\eta_j$, hence $\hat{\pi}_0^H$ and  $\hat{\pi}_0^G$ are identical for a fixed set of tuning parameters. Otherwise for each $j\in \mathcal{J}$, $\eta_j\geq \tau_j$ and thus $\hat{\pi}_0^H$ dominates $\hat{\pi}_0^G$ almost surely. This small amount of conservative bias enables us to prove the desired result theoretically without the impractical assumption taken up in Chen and Doerge (2020). However, we already mentioned that $\hat{\pi}_0^G$ and $\hat{\pi}_0^S$ are identical under the same assumption and thus a good point of investigation is to validate whether the conservative bias present in $\hat{\pi}_0^S$ is more than the bias present in $\hat{\pi}_0^H$.

\subsection{Bias of the estimator}
From the definition of $\hat{\pi}_0^H$, Bias$(\hat{\pi}_0^H)=$ $(1/n)$Bias$[\beta(\tau_j)]$. Now using the same notations as in Section 2

\begin{equation}\label{eq:Ebeta}
    E[\beta(\tau_j)]=\frac{1}{m(1-\eta_j)}+\frac{1-\tau_j}{m(1-\eta_j)}\left[\sum_{i\in \mathcal{T}}\frac{1-F_i(\lambda_{ij})}{1-\lambda_{ij}} +\sum_{{i\in\mathcal{F}}} \frac{1-G_i(\lambda_{ij})}{1-\lambda_{ij}}\right].
\end{equation}
As defined in Algorithm 1, $\lambda_{ij}\in S_i$ and hence $F_i(\lambda_{ij})=\lambda_{ij}$ for each $i\in \mathcal{I}$ and $j\in \mathcal{J}$. Thus $\sum_{i\in \mathcal{T}} (1-F_i(\lambda_{ij}))/(1-\lambda_{ij})$ in \ref{eq:Ebeta} is equal to $m_0$ and $E(m_0)=m(1-\pi_0)$ as discussed in Section 2. Therefore

\begin{equation}\label{eq:Biasbeta}
    Bias[\beta(\tau_j)]=\frac{1}{m(1-\eta_j)}+\left(1-\frac{1-\tau_{j}}{1-\eta_j}\right)\pi_0+\frac{1-\tau_j}{m(1-\eta_j)}\sum_{i\in \mathcal{F}}\frac{1-G_i(\lambda_{ij})}{1-\lambda_{ij}}.
\end{equation}
With $\tau=\tau_j$ in $B_2$, the bias structure of $\hat{\pi}_0^S$ we get 

\begin{equation}\label{eq:Biasstorey}
    Bias[\hat{\pi}_0(\tau_j)]=\frac{1}{m(1-\tau_j)}+\frac{1}{m(1-\tau_j)}\sum_{i\in \mathcal{T}}[\tau-F_i(\tau_j)] +\frac{1}{m(1-\tau_j)}\sum_{i\in \mathcal{F}}[1-G_i(\tau_j)].
\end{equation}
Note that under that assumption $\tau_j=\eta_j$ for each $j\in \mathcal{J}$, the bias in $\beta(\tau_j)$ is identical to $\hat{\pi}_0^j$ as expected. For the general setting, the first term in (\ref{eq:Biasbeta}) dominates the first term in (\ref{eq:Biasstorey}). The other terms are not directly comparable. Numerical comparison of bias of the new estimator with the estimator with $\hat{\pi}_0^S(0.5)$ under different simulation settings will provide a deeper insight. 

\section{FDR control by adaptive BH algorithm}
Denote the set of available $p$-values $\{p_1, p_2, ..., p_m\}$ by $\bm{p}$. Let us introduce the following notation.

\begin{equation*}
    \bm{p}_{k}=\{p_1, ..., p_{k-1}, 0, p_{k+1}, ..., p_{m}\}
\end{equation*}
Accordingly we define $\beta_k(\tau_j)$ and $\hat{\pi}_{0k}^H$ as $\beta(\tau_j)$ and $\hat{\pi}_{0}^H$ computed by replacing $\bm{p}$ by $\bm{p}_k$, respectively for each $k\in \mathcal{I}$ and $j\in \mathcal{J}$. \\

We briefly revisit the BH algorithm. First $\bm{p}$ is ordered $p_{(1)}\leq p_{(2)}\leq ... \leq p_{(m)}$ and the hypothesis corresponding to $p_{(i)}$ is renamed $H_{(i)}$, $i\in \mathcal{I}$. Define $\hat{k}=\max_{1\leq i\leq m}\{i:p_{(i)}\leq i\alpha/m\}$ for some prefixed level $\alpha\in (0,1)$. If such a $\hat{k}$ exists, we reject $H_{(1)}, H_{(2)}, ..., H_{(\hat{k})}$ and accept the remaining hypotheses. If such a $\hat{k}$ does not exist, then we do not reject any hypotheses. The adaptive version of BH algorithm (Henceforth ABH) takes the choice of $\hat{k}=\max_{1\leq i\leq m}\{i:\hat{\pi}_0^H p_{(i)}\leq i\alpha/m\}$ and proceeds in similar line with the BH algorithm. The adjusted $p$-values of the BH algorithm and the BHH algorithm are discussed in Heyse (2011). The new estimator for $\pi_0$ can similarly be plugged into the BHH algorithm to construct adaptive BHH algorithm (Henceforth ABHH) (see Chen et al. 2018). However, the BHH and hence ABHH algorithms are difficult to express in terms of linear cut-offs for ordered $p$-values as BH and BHH. Thus the technique we opt for proving FDR control by ABH cannot be directly applied to ABHH. Due to the improved performance of BHH over BH algorithm, numerical investigation regarding FDR control and real life application of ABHH are worth performing.\\   

\noindent \textbf{Theorem 1}: If the $p$-values are independent of each other and BH procedure controls the FDR at a prefixed level $\alpha\in (0,1)$, then the $\hat{\pi}_0^H$ plugged-in ABH procedure controls the FDR at $\alpha$. \\

\noindent \textbf{Proof of Theorem 1}: Obviously $\hat{\pi}_0^H=H(p_1, p_2, ..., p_m)$, a function of the available $p$-values for a fixed set of tuning parameters. The functional reciprocal of $\hat{\pi}_0^H$ be $G(p_1, p_2, ..., p_m)$, for some function $G: (0,1)^m\to R^+$. Noting that the BH procedure is a linear step-up procedure, we apply Theorem 11 of Blanchard and Roquain (2009) (A particular case was done earlier in Benjamini et al. (2006)). Thus to prove Theorem 1 we need to validate the following sufficient conditions on $\hat{\pi}_0^H$.
\begin{enumerate}
    \item $G$ is a component-wise non-decreasing function.
    \item For each $k\in \mathcal{T}$, $E\left[1/\hat{\pi}_{0k}^H\right]\leq1/\pi_0$.
\end{enumerate}
For $p_i^\prime\geq p_i$ for some $i\in \mathcal{I}$, $I(p_i^\prime>\lambda_{ij})\geq I(p_i\geq \lambda_{ij})$ for each $j\in \mathcal{J}$. Other terms in the expression of $\hat{\pi}_0^H$ are free of any $p_i$. From the above facts it is trivial to justify the first condition.\\

The second condition is the crucial one. From the definition of $\hat{\pi}_0^H$ in Algorithm 1 and the notations introduced earlier in this section, using AM-HM inequality we get

\begin{equation*}
    \frac{1}{\hat{\pi}_{0k}^H}=\frac{n}{\sum_{j\in \mathcal{J}}\beta_k(\tau_j)}\leq \frac{1}{n}\sum_{j\in \mathcal{J}}\frac{1}{\beta_k(\tau_j)}
\end{equation*}
and hence

\begin{equation*}
    E\left[\frac{1}{\hat{\pi}_{0k}^H}\right]\leq \frac{1}{n}\sum_{j\in \mathcal{J}}E\left[\frac{1}{\beta_k(\tau_j)}\right].
\end{equation*}
Thus for proving the second condition, it is sufficient to show for each $j\in \mathcal{J}$ that

\begin{equation}\label{eq:secondcondition}
    E\left[\frac{1}{\beta_k(\tau_j)}\right]\leq E\left[\frac{1}{\pi_0}\right].
\end{equation}
Now for a fixed $j\in\mathcal{J}$

\begin{equation*}
    \beta(\tau_j)=\frac{1}{m(1-\eta_j)}+\frac{1-\tau_j}{m(1-\eta_j)}\left[\sum_{i \in \mathcal{I}}\frac{I(p_i>\lambda_{ij})}{1-\lambda_{ij}}\right]
\end{equation*}
and hence 

\begin{equation}\label{eq:betak}
    \beta_k(\tau_j)=\frac{1}{m(1-\eta_j)}+\frac{1-\tau_j}{m(1-\eta_j)}\left[\sum_{i \in \mathcal{I}_k}\frac{I(p_i>\lambda_{ij})}{1-\lambda_{ij}}\right].
\end{equation}
Here $\mathcal{I}_k=\{1, ..., k-1, 0, k+1, ..., m\}$. For $k\in\mathcal{T}$, denote $\mathcal{T}-\{k\}$ by $\mathcal{T}_k$. Replacing $I(p_i>\lambda_{ij})/(1-\lambda_{ij})$ for $i\in \mathcal{I}_k-\mathcal{T}_k$ by 0 in equation (\ref{eq:betak}) we get

\begin{equation}\label{eq:mbetak}
    m\beta_k(\tau_j)\geq \frac{1}{1-\eta_j}+\frac{1-\tau_j}{1-\eta_j}\sum_{i\in \mathcal{T}_k}\frac{I(p_i>\lambda_{ij})}{1-\lambda_{ij}}
\end{equation}
and hence
\begin{equation}\label{eq:Eonebybeta}
    E\left[\frac{1}{\beta_k(\tau_j)}\right]\leq m(1-\eta_j)\, E\left[\frac{1}{1+(1-\tau_j)\sum_{i\in \mathcal{I}}\frac{I(p_i>\lambda_{ij})}{1-\lambda_{ij}}}\right].
\end{equation}
For all $i\in\mathcal{I}$, $\lambda_{ij}\geq\tau_j$ and thus $(1-\lambda_{ij})^{-1}\geq (1-\tau_j)^{-1}$. Replacing $(1-\lambda_{ij})^{-1}$ by $(1-\tau_j)^{-1}$ in right hand side (RHS) of 
equation (\ref{eq:Eonebybeta}) we get

\begin{equation}
    E\left[\frac{1}{1+(1-\tau_j)\sum_{i\in \mathcal{I}}\frac{I(p_i>\lambda_{ij})}{1-\lambda_{ij}}}\right]\leq E\left[\frac{1}{1+\sum_{i\in \mathcal{T}_k}I(p_i>\lambda_{ij})}\right]
\end{equation}
and hence

\begin{equation}\label{eq:Eonebybeta2}
    E\left[\frac{1}{\beta_k(\tau_j)}\right]\leq m(1-\eta_j)\,E\left[\frac{1}{1+\sum_{i\in \mathcal{T}_k}I(p_i>\lambda_{ij})}\right].
\end{equation}
Since $\eta_j=\max_{i\in \mathcal{I}}\lambda_{ij}$, inequality in equation (\ref{eq:Eonebybeta2}) is preserved by replacing $\lambda_{ij}$ by $\eta_j$ for each $i\in \mathcal{I}$ RHS of (\ref{eq:Eonebybeta2}). 

\begin{equation}\label{eq:Eonebybeta3}
    E\left[\frac{1}{\beta_k(\tau_j)}\right]\leq m(1-\eta_j)\,E\left[\frac{1}{1+\sum_{i\in \mathcal{T}_k}I(p_i>\eta_j)}\right]
\end{equation}
For each $i\in \mathcal{I}$, $p_i\geq U_i\sim$ Uniform$(0,1)$ almost surely. Therefore $I(p_i>\eta_j)\geq I(U_i\geq \eta_j)$ almost surely. Since $p_i$'s are independent of each other $U_i$'s are considered to be so. Thus 

\begin{equation}\label{eq:Eonebybeta4}
   E\left[\frac{1}{\beta_k(\tau_j)}\right]\leq m(1-\eta_j)\,E\left[\frac{1}{1+\sum_{i\in \mathcal{T}_k}I(U_i>\eta_j)}\right]. 
\end{equation}
Since $U_i$'s are independent, $B=\sum_{i\in \mathcal{T}_k}I(U_i>\eta_j)\sim$ Binomial$(m_0-1,1-\eta_j)$. Applying Lemma 1 of Benjamini et al. (2006) in equation (\ref{eq:Eonebybeta4}) we get 

\begin{equation}\label{eq:ineq}
    E\left[\frac{1}{1+B}\right]\leq\frac{1}{m_0(1-\eta_j)}.
\end{equation}
Thus using (\ref{eq:ineq}) in equation (\ref{eq:Eonebybeta4}) we get the result in (\ref{eq:secondcondition}). Hence the proof.

\section{Remarks}
This article contains some theoretical developments regarding the proposed estimator for the proportion of true null hypotheses. Work on a data-driven choice of the tuning parameters is currently ongoing. Upon arriving at a suitable selection technique, future work consists of performance exploration through extensive simulation studies. The next version will also accommodate a real life application to establish the new methodology. We are expecting to report the final results in a revised version of the current article within a short period of time.  

\section*{References}
\begin{small}
Benjamini, Y., and Hochberg, Y. (1995), ``Controlling the false discovery rate: a practical and powerful approach to multiple testing," \textit{Journal of the Royal Statistical Society: Series B (Statistical Methodology)}, 289-300.\\
\\
Benjamini, Y., Krieger, A. M., and Yekutieli, D. (2006), ``Adaptive linear step-up procedures that control the false discovery rate," \textit{Biometrika}, 93(3), 491-507.\\
\\
Biswas, A. (2020), ``Estimating the proportion of true null hypotheses based on sum of $p$-values and application in microarray data," \textit{Communications in Statistics - Simulation and Computation}, DOI: 10.1080/03610918.2020.1800036\\
\\
Biswas, A. (2020), ``Regarding Paper ``Multiple testing with discrete data: Proportion of true null hypotheses and two adaptive FDR procedures" by Xiongzhi Chen, Rebecca W. Doerge, and Joseph F. Heyse," \textit{Biometrical journal}. Biometrische Zeitschrift.\\
\\
Blanchard, G., and Roquain, E. (2009), ``Adaptive false discovery rate control under independence and dependence," \textit{Journal of Machine Learning Research}, 10, 2837–2871.\\
\\
Chen, X., and Doerge, R. (2014), ``Generalized estimators for multiple testing: Proportion of true nulls and false discovery rate," Retrieved from \href{http://arxiv.org/abs/1410.4274}{http://arxiv.org/abs/1410.4274}.\\
\\
Chen, X., Doerge, R. W., and Heyse, J. F. (2018), ``Multiple testing with discrete data: Proportion of true null hypotheses and two adaptive FDR procedures," \textit{Biometrical Journal}, 60(4), 761-779.\\
\\
Chen, X. (2020), ``False discovery rate control for multiple testing based on discrete p‐values," \textit{Biometrical Journal}, 62(4), 1060-1079.\\
\\
Chen, X., and Doerge, R. W. (2020), ``Comments on Dr. Aniket Biswas' Letter to the Editor," \textit{Biometrical Journal}.\\
\\
Cheng, Y., Gao, D., and Tong, T. (2015), ``Bias and variance reduction in estimating the proportion of true-null hypotheses," \textit{Biostatistics}, 16(1), 189-204.\\
\\
Dialsingh, I., Austin, S. R., and Altman, N. S. (2015), ``Estimating the proportion of true null hypotheses when the statistics are discrete,"\textit{Bioinformatics}, 31(14), 2303-2309.\\
\\
Heyse, J. F. (2011), ``A false discovery rate procedure for categorical data," In M. Bhattacharjee, S. K. Dhar, \& S. Subramanian (Eds.), \textit{Recent advancesin biostatistics: False discovery rates, survival analysis, and related topics} (Vol. 4, chap. 3, pp. 43–58). Singapore, SI: World Scientific.\\
\\
Jiang, H., and Doerge, R. W. (2008), "Estimating the proportion of true null hypotheses for multiple comparisons," \textit{Cancer Informatics}, 6, 117693510800600001.\\
\\
Langaas, M., Lindqvist, B. H., and Ferkingstad, E. (2005), ``Estimating the proportion of true null hypotheses, with application to DNA microarray data,"\textit{ Journal of the Royal Statistical Society: Series B (Statistical Methodology), }67(4), 555-572.\\
\\
Sarkar, S. K. (2008), ``On methods controlling the false discovery rate," \textit{Sankhyā: The Indian Journal of Statistics, Series A} (2008-), 135-168.\\
\\
Storey, J. D. (2002), ``A direct approach to false discovery rates," \textit{Journal of the Royal Statistical Society: Series B (Statistical Methodology), }64(3), 479-498.\\
\\
Storey, J. D., Taylor, J. E., and Siegmund, D. (2004), ``Strong control, conservative point estimation and simultaneous conservative consistency of false discovery rates: a unified approach," \textit{Journal of the Royal Statistical Society: Series B (Statistical Methodology)}, 66(1), 187-205.\\
\\
Wang, H. Q., Tuominen, L. K., and Tsai, C. J. (2010), ``SLIM: a sliding linear model for estimating the proportion of true null hypotheses in datasets with dependence structures,"\textit{ Bioinformatics}, 27(2), 225-231.\\
\end{small}
\end{document}